\documentclass[11pt]{article}

\usepackage{amsthm}%
\usepackage{amssymb}%
\usepackage{amsmath}%

\usepackage{indentfirst}%
\newtheoremstyle{cute}
    {} 
    {} 
    {\slshape} 
    {} 
    {\bfseries} 
    { } 
    { } 
    {} 

\newtheoremstyle{definitioncute}
    {} 
    {} 
    {\upshape} 
    {} 
    {\bfseries} 
    { } 
    { } 
    {} 

\theoremstyle{cute}
\newtheorem{thm}{Theorem}

\newtheorem*{namedtheorem}{\theoremname}
\newcommand{\theoremname}{Theorem}
    {\end{namedtheorem}}

\theoremstyle{definitioncute}
\theoremstyle{remark}

\newcounter{chislo}
    {\begin{list}{(\roman{chislo})}{\usecounter{chislo}%
        \setlength{\topsep}{.5ex}\setlength{\labelwidth}{50pt}%
        \setlength{\parsep}{0pt}\setlength{\itemsep}{0pt}%
        \setlength{\leftmargin}{30pt}}}%
    {\end{list}}
\makeatletter

\renewcommand{\@seccntformat}[1]{\csname the#1\endcsname.\hspace{.5em}}

\renewcommand{\section}{\@startsection
    {section}
    {1}
    {0em}
    {1.5\baselineskip plus .5\baselineskip minus .5\baselineskip}
    {.5\baselineskip plus .2\baselineskip minus 0\baselineskip}
    {\normalfont\normalsize\scshape\large\centering}}

\renewcommand{\subsection}{\@startsection
    {subsection}
    {2}
    {0em}
    {.5\baselineskip plus .5\baselineskip minus 0\baselineskip}
    {.5\baselineskip}
    {\normalfont\normalsize\scshape}}

\makeatother

\renewcommand{\title}[1]{\begin{center}\bfseries\Large #1\end{center}\smallskip}
\renewcommand{\author}[1]{\begin{center}\scshape\normalsize #1\end{center}}
\begin{document}

\title{Boundedly Simple Groups Have Trivial Bounded Cohomology}

\author{Igor V.\ Erovenko}

\begin{center}
May 2, 2004
\end{center}

The goal of this short note is to observe that the singular part of the second
bounded cohomology group of boundedly simple groups constructed in
\cite{muranov} is trivial. Recall that a group $G$ is called
\emph{$m$-boundedly simple} if every element of $G$ can be represented as a
product of at most $m$ conjugates of $g$ or $g^{-1}$ for any $g \in G$.

We recall that bounded cohomology $H^*_b(G)$ of a group $G$ (we will be
considering only cohomology with coefficients in the additive group of reals
$\mathbb R$ with trivial action, so in our notations for cohomology the
coefficient module will be omitted) is defined using the complex
\[
\cdots \longleftarrow C^{n+1}_b(G) \stackrel{\delta^{n}_b}{\longleftarrow}
C^{n}_b(G) \longleftarrow \cdots \longleftarrow C^{2}_b(G)
\stackrel{\delta^{1}_b}{\longleftarrow} C^{1}_b(G)
\stackrel{\delta^{0}_b=0}{\longleftarrow} \mathbb{R}
\stackrel{\delta^{-1}_b=0}{\longleftarrow} 0
\]
of bounded cochains $f \colon G \times \cdots \times G \to \mathbb R$, and
$\delta^{n}_b = \delta^{n} |_{C^{n}_b(G)}$ is the bounded differential
operator. Since $H_b^0 (G) = \mathbb R$ and $H^1_b(G) = 0$ for any group $G$,
investigation of bounded cohomology starts in dimension 2. One observes that
$H_b^2(G)$ contains a subspace $H_{b,2}^2(G)$ (called the \emph{singular part}
of the second bounded cohomology group), which has a simple algebraic
description in terms of quasicharacters and pseudocharacters, and the quotient
space $H^2_b(G)/H_{b,2}^2(G)$ is canonically isomorphic to the bounded part of
the ordinary cohomology group $H^2(G)$. See \cite{grigorchuk-95} for background
and available results on bounded cohomology of groups.

A function $F \colon G \to \mathbb{R}$ is called a \emph{quasicharacter} if
there exists a constant $C_F \geqslant 0$ such that
\[
| F(xy) - F(x) - F(y) | \leqslant C_F \quad \text{for all} \ x, y \in G.
\]
A function $f \colon G \to \mathbb{R}$ is called a \emph{pseudocharacter} if
$f$ is a quasicharacter and in addition
\[
f(g^n)=nf(g) \quad \text{for all } g \in G \text{ and } n \in \mathbb{Z}.
\]
We use the following notation: $X(G)=$ the space of additive characters $G \to
\mathbb{R}$; $QX(G)=$ the space of quasicharacters; $PX(G)=$ the space of
pseudocharacters; $B(G)=$ the space of bounded functions. Then
\begin{equation}\label{eq:pseudo}
H^{2}_{b,2}(G) \cong QX(G)/ (X(G) \oplus B(G)) \cong PX(G)/ X(G)
\end{equation}
as vector spaces (cf.\ \cite[Proposition~3.2 and Theorem~3.5]{grigorchuk-95}).
Special interest in $H^2_{b,2}$ is motivated in part by its connections with
other structural properties of groups such as commutator length \cite{bavard}
and bounded generation \cite{grigorchuk-95}.

\begin{thm}
If $G$ is a boundedly simple group, then $H^2_{b,2}(G) = 0$.
\end{thm}
\begin{proof}
In view of \eqref{eq:pseudo} it suffices to show that the group $G$ does not
have any nontrivial pseudocharacters. First, we observe that every
pseudocharacter is constant on conjugacy classes. Indeed, suppose that $f \in
PX(G)$ and $| f(g x g^{-1}) - f(x) | = a > 0$ for some $x$, $g \in G$. Then on
the one hand
\[
| f(g x^n g^{-1}) - f(x^n) | = | f(g x^n g^{-1}) - f(x^n) - f(g) - f(g^{-1}) |
\leqslant 2 C_f
\]
is bounded independent of $n$, on the other hand
\[
| f(g x^n g^{-1}) - f(x^n) | = n | f(g x g^{-1}) - f(x) | = na \to \infty \quad
\text{as } n \to \infty,
\]
whence a contradiction.

Suppose that $G$ is $m$-boundedly simple. Then every element $x$ of $G$ can be
written in the form
\[
x = g_1 \cdots g_k
\]
where $k \leqslant m$ and every $g_i$ is a conjugate of either $g$ or $g^{-1}$
for some fixed $g \in G$, whence $|f(g_i)| = |f(g)|$ for all $i = 1, \dots ,
k$. Then
\begin{eqnarray*}
|f(x)| & = & |f(g_1 \cdots g_k) - f(g_1) - \cdots - f(g_k) + f(g_1) + \cdots +
f(g_k)| \\
& \leqslant & |f(g_1 \cdots g_k) - f(g_1) - \cdots - f(g_k)| + |f(g_1)| +
\cdots + |f(g_k)| \\
& \leqslant & (m-1) C_f + m |f(g)|
\end{eqnarray*}
which implies that $f$ is bounded on $G$, hence must be trivial.
\end{proof}

\begin{flushleft}
Department of Mathematical Sciences\\
University of North Carolina at Greensboro\\
Greensboro NC 27402\\
E-mail: igor@uncg.edu
\end{flushleft}

\end{document}